\documentclass[11pt,reqno,tbtags]{amsart}
\usepackage{amssymb}
\numberwithin{equation}{section}

\newtheorem{theorem}{Theorem}[section]
\newtheorem{lemma}[theorem]{Lemma}

\newtheorem{condition}[theorem]{Condition}

\theoremstyle{definition}

\newtheorem{remark}[theorem]{Remark}
\newtheorem{remarks}[theorem]{Remarks}

\newtheorem*{acks}{Acknowledgements}

\theoremstyle{remark}

\newenvironment{romenumerate}{\begin{enumerate}
 }{\end{enumerate}}

\newcounter{thmenumerate}
\newenvironment{thmenumerate}
{\setcounter{thmenumerate}{0}%
 \def\item{\ifnum\thethmenumerate=0\else\par\fi 
 \addtocounter{thmenumerate}{1}\textup{(\roman{thmenumerate})\enspace}}}
{}

\newcounter{xenumerate}

\newcommand{\refT}[1]{Theorem~\ref{#1}}
\newcommand{\refC}[1]{Corollary~\ref{#1}}
\newcommand{\refCN}[1]{Condition~\ref{#1}}
\newcommand{\refL}[1]{Lemma~\ref{#1}}

\newcommand{\refS}[1]{Section~\ref{#1}}

\newcommand{\refand}[2]{\ref{#1} and~\ref{#2}}

\begingroup
  \divide\count255 by 60
  \multiply\count255 by -60
  \advance\count255 by \time
  \ifnum \count255 < 10 \xdef\klockan{\the\count1.0\the\count255}
  \else\xdef\klockan{\the\count1.\the\count255}\fi
\endgroup

\newcommand\nopf{\qed}   

\newcommand{\sumr}{\sum_{r=0}^\infty}
\newcommand{\sumoo}{\sum^\infty}
\newcommand{\sumin}{\sum_{i=1}^n}

\newcommand\set[1]{\ensuremath{\{#1\}}}

\newcommand\bigpar[1]{\bigl(#1\bigr)}
\newcommand\Bigpar[1]{\Bigl(#1\Bigr)}

\def\rompar(#1){\textup(#1\textup)}    

\newcommand\parfrac[2]{\Bigpar{\frac{#1}{#2}}}

\def\xexp(#1){e^{#1}}
\newcommand\ceil[1]{\lceil#1\rceil}

\newcommand\ntoo{\ensuremath{{n\to\infty}}}

\newcommand\xtoo{\ensuremath{{x\to\infty}}}

\newcommand\iid{i.i.d.\spacefactor=1000}    
\newcommand\ie{i.e.\spacefactor=1000}
\newcommand\eg{e.g.\spacefactor=1000}

\newcommand\cf{cf.\spacefactor=1000}
\newcommand{\as}{a.s.\spacefactor=1000}

\newcommand{\tend}{\longrightarrow}

\newcommand\pto{\overset{\mathrm{p}}{\tend}}

\newcommand\Op{O_{\mathrm p}}
\newcommand\op{o_{\mathrm p}}

\newcommand\bbN{\mathbb N}

\newcommand\bbZ{\mathbb Z}

\newcommand\E{\operatorname{\mathbb E{}}}
\renewcommand\P{\operatorname{\mathbb P{}}}

\newcommand\Exp{\operatorname{Exp}}
\newcommand\Po{\operatorname{Po}}
\newcommand\Bi{\operatorname{Bi}}

\newcommand\ga{\alpha}

\newcommand\gd{\delta}

\newcommand\gl{\lambda}

\newcommand\eps{\varepsilon}

\renewcommand\phi{\varphi}

\newcommand\ett[1]{\boldsymbol1[#1]} 

\def\[#1]{[\![#1]\!]}

\renewcommand{\=}{:=}

\newcommand\rhs{right hand side}

\newcommand\ltx{<\!}
\newcommand\gex{\ge\!}
\newcommand\kcore{$k$-core}
\newcommand\gnp{\ensuremath{G(n,p)}}
\newcommand\gnm{\ensuremath{G(n,m)}}
\newcommand\gnl{\ensuremath{G(n,\gl/n)}}
\newcommand\gnd{\ensuremath{G(n,(d_i)_1^n)}}
\newcommand\gndx{\ensuremath{G^*(n,(d_i)_1^n)}}
\newcommand\crk{\ensuremath{\mathrm{Core}_k}}
\newcommand\crkx{\ensuremath{\mathrm{Core}^*_k}}
\newcommand\nn{^{(n)}}
\newcommand\et{e^{-t}}
\newcommand\tij{T_i}
\newcommand\psij{\psi_{j}}
\newcommand\psik{\psi_{k}}
\newcommand\psiki{\psi_{k-1}}

\newcommand\glk{\lambda_{k}}
\newcommand\muxx[1]{\mu^{#1}}
\newcommand\pmax{\widehat p}
\newcommand\tnxt{N^{(\ceil x)}}

\newcommand\whp{\textbf{whp}}

\begin{document}
\title
{A simple solution to the $k$-core problem}

\date{24 August 2005 (\today)}

\author{Svante Janson}
\address{Department of Mathematics, Uppsala University, PO Box 480,
SE-751~06 Uppsala, Sweden}
\email{svante.janson@math.uu.se}
\urladdr{http://www.math.uu.se/\~{}svante/}

\author{Malwina J. Luczak}
\address{Department of Mathematics, London School of Economics,
  Houghton Street, London WC2A 2AE, United Kingdom}
\email{malwina@planck.lse.ac.uk}
\urladdr{http://www.lse.ac.uk/people/m.j.luczak@lse.ac.uk/}


\begin{abstract} 
We study 
the $k$-core of a random (multi)graph on $n$
vertices with a given degree sequence.
We let \ntoo. Then, under some regularity conditions on the degree
sequences,
we give conditions on the asymptotic shape of the degree sequence
that imply that
with high probability the $k$-core is empty, 
and other conditions that imply that with high probability
the $k$-core is non-empty and the sizes of its vertex and edge sets
satisfy a law of large numbers; 
under suitable assumptions these are the only two possibilities.
In particular, we recover the result by Pittel, Spencer and Wormald
\cite{psw96} on the existence and size of a $k$-core in $G(n,p)$ and
$G(n,m)$,
see also Molloy~\cite{Molloy05} and Cooper~\cite{c04}.

Our method is based on the properties
of empirical distributions of independent random variables, and leads
to simple proofs.
\end{abstract}

\maketitle

\section{Introduction}\label{S:intro}

Let $k \ge 2$ be a fixed integer. The
$k$-core of a graph $G$ 
is the largest induced
subgraph of $G$ with minimum vertex degree at least $k$. 
The question whether a non-empty $k$-core exists in a random graph has
attracted a lot of attention over the past fifteen years. There have 
by now been quite a number of studies for the Bernoulli random graph $\gnp$
with $n$ vertices and edge
probability $p$, and for the uniformly random graph $\gnm$ with $n$
vertices and $m$ edges (see~\cite{c04,l91,Molloy05,psw96} and references
therein). Recently, Fernholz and Ramachandran~\cite{fr03,fr04} have
considered the $k$-core of a random graph with a specified degree
sequence. More generally, Cooper~\cite{c04} studies cores of
random uniform hypergraphs with a given degree sequence.   
Yet more generally, Molloy~\cite{Molloy05} considers cores in random
structures such as the uniform hypergraph and satisfiability of
boolean formulas (see also references therein).

For a constant $\mu > 0$, let $\Po (\mu)$ denote a Poisson random
variable with mean $\mu$. Given $\mu > 0$ and $j \in
\bbZ^+$, let $\psij(\mu)\=\P\bigpar{\Po(\mu)\ge j}$. Also, let
$\glk\=\min_{\mu>0}\mu/\psiki(\mu)$; and for $\gl>\glk$, we use
$\mu_k (\gl) > 0$ to denote the largest solution to $\mu/\psiki(\mu)=\gl$.

In~\cite{psw96}, Pittel, Spencer and Wormald discovered that for $k
\ge 3$, $\gl = \glk$ 
is the threshold for the appearance of a nonempty $k$-core in the graph
\gnl{} (or, equivalently, $m = n\glk/2$ is the threshold in the graph $\gnm$).
Their strategy was to analyse an edge deletion algorithm that finds
the $k$-core in a graph, showing that the corresponding random process
is well approximated by the solution to a system of differential equations. 
The proof is rather long and complicated, 
and involves counting formulae for the
number of graphs with a given degree sequence.
For an analysis that uses a slightly modified version of their
deletion algorithm and differs in some other important technical
details too, see ~\cite{l05}.

Fernholz and Ramachandran  \cite{fr03,fr04} use different techniques to
study the existence of a large $k$-core in a
random graph with a given degree sequence. Their core-finding  
algorithm is basically identical to ours, but they analyse it in
quite a different way; they also compare their result
to a corresponding result for branching processes.

Cooper~\cite{c04} has studied the $k$-core of a uniform
multihypergraph with a given degree sequence. His method involves
analysing a constructive algorithm generating the
multihypergraph and its core, and inductively applying Azuma's
inequality over time periods of length $n^{2/3} \Delta^{4/3} \log n$,
where $\Delta$ is the initial maximum degree. 

Molloy~\cite{Molloy05} gave another proof of the sharp threshold
for the $k$-core, analysing a multi-round vertex and edge deletion
algorithm via a branching process type argument.

Kim~\cite{kim} considers cores in a ``Poisson cloning'' model
of a random graph, which is somewhat different from \gnp. The
slides~\cite{kim} present a sketch argument, without precise error
bounds, showing that the critical threshold for the emergence of a $k$-core
agrees with the threshold in \gnp.

Darling and Norris~\cite{dn05} analyse
cores in a different, weighted, Poisson model of a random
hypergraph. Their method involves establishing a differential equation
approximation for the Markov chain representing a suitable
deletion algorithm. The threshold for \gnp{}
follows as a corollary to their main result.

Also see Cain and Wormald \cite{CW}, who use differential equations
to analyse the $k$-core threshold and the properties of the degree
sequence of the giant $k$-core in a different model of a random
graph. They make corresponding statements for \gnm\ as a corollary.

In this paper, we present a simple solution to the $k$-core
problem. Unlike \cite{dn05} and~\cite{psw96}, we do not use
differential equations, but rely solely 
on the convergence of empirical distributions of independent
random variables. Apart from \gnp{} and \gnm, we are also able to handle the
uniformly random graph with a given degree sequence under some
regularity conditions 
similar to \cite{c04,fr03,fr04}. 
In contrast to~\cite{dn05,psw96} we
do not require counting formulae for graphs but, like~\cite{c04}
and~\cite{fr03,fr04}, work
directly in the configuration model used to construct the random
graph, exposing the edges one by one as they are needed.

We shall now state the result  concerning the
emergence of the $k$-core in the random graphs \gnp{} and \gnm.
Given a graph $G$, let 
$v(G)$ and $e(G)$ denote the sizes of the vertex and edge sets of $G$
respectively. 
We consider asymptotics as \ntoo, and 
say that an event holds \whp{} (\emph{with high probability}), 
if it holds with probability tending to 1 as 
$n\to\infty$.

We shall use $\Op$ and $\op$ in the standard way (see \eg{}
Janson, {\L}uczak and Ruci\'nski~\cite{JLR}); for example,
if $(X_n)$ is a sequence of random variables, then
$X_n=\Op(1)$ means ``$X_n$ is bounded in probability'' and
$X_n=\op(1)$ means that $X_n \pto0$.

\begin{theorem}[Pittel, Spencer and Wormald \cite{psw96}]
  \label{T1}\hskip 0pt plus 3pt
Consider the random  graph \gnl, where $\gl>0$ is fixed.
Let $k\ge2$ be fixed and let $\crk = \crk (n,\gl)$ be the \kcore{} of \gnl.
  \begin{romenumerate}
\item
If $\gl<\glk$ and $k\ge3$,  then 
$\crk$ is empty \whp.
\item
If $\gl>\glk$,
then
\whp{} \crk{} is non-empty, and 
$v(\crk)/n\pto \psik(\mu_k(\gl))$, $e(\crk)/n\pto
\mu_k(\gl)\psiki(\mu_k(\gl))/2=\mu_k(\gl)^2/(2\gl)$.
  \end{romenumerate}

The same results hold for the random graph \gnm, for any sequence
$m=m(n)$ with $2m/n\to\gl$.
\end{theorem}

Part (i) does not hold for $k=2$. Here $\lambda_2=1$ and for $0<\gl<1$
there is a positive limiting probability that there are cycles
(as shown already by Erd\H os and R\'enyi \cite{ER60}), 
and thus a non-empty 2-core. Nevertheless, in this
case $e(\crk)=\Op(1)$ and $v(\crk)=\Op(1)$, so the core is small;
\cf{} \refT{TD}(i) below.

\begin{acks}
  This research was mainly done during a visit by MJL to Uppsala
  University in April 2005, sponsored by the LSE Nordic Exchange
  Scheme.
\end{acks}

\section{Multigraphs}\label{Smulti}

It will be convenient to work with \emph{multigraphs}, that is to allow
multiple edges and loops.
In particular, we shall use the following type of random multigraph.

Let $n \in \bbN$ and let $(d_i)_1^n=(d_i\nn)_1^n$ be a sequence of
non-negative integers such that  $\sumin d_i$ is even.
We define a \emph{random multigraph with given degree
sequence $(d_i)_1^n$}, denoted by \gndx, 
by the configuration model (see \eg{} \cite{b01}): 
take a set of $d_i$ half-edges for each vertex 
$i$, and combine the half-edges into pairs by a uniformly random
matching of the set of all half-edges. 
Note that \gndx{} does not have exactly the uniform
distribution over all multigraphs with the given degree sequence;
there is a weight with a factor $1/j!$ for every edge of multiplicity $j$, and 
a factor $1/2$ for every loop, see \cite[\S1]{giant}. However,
conditioned on the multigraph being a (simple) graph, we obtain a
uniformly distributed random graph with the given degree sequence,
which we denote by \gnd.

\begin{remark}
The distribution of \gndx{} is the same as the one obtained by sampling the
edges as ordered pairs of vertices uniformly with replacement, and
then conditioning on the vertex degrees being correct.  
\end{remark}

Let us write $2m\=\sumin d_i$, so that $m=m(n)$ is the number of
edges in the multigraph \gndx.
We will let \ntoo, and assume that we are given
$(d_i)_1^n$ satisfying the following regularity conditions,
\cf{} Molloy and Reed \cite{MR95}.
\begin{condition}\label{C1}
For each $n$,  $(d_i)_1^n=(d_i\nn)_1^n$ is a sequence of non-negative
integers such that $\sumin d_i$ is even and, for some probability distribution
$(p_r)_{r=0}^\infty$ independent of $n$, 
\begin{romenumerate}
  \item
$\#\set{i:d_i=r}/n\to p_r$ for every $r\ge0$ as \ntoo;
\item 
$\gl\=\sum_r r p_r\in(0,\infty)$;
\item
$2m/n\to\gl$ as \ntoo.
\end{romenumerate}
\end{condition}

We shall consider thinnings of the vertex degrees in \gndx.
Let $W$ be a random variable with the distribution $\P(W=r)=p_r$.
(This is the asymptotic distribution of the vertex degrees in \gndx.)
For $0\le p\le1$ we let $W_p$ be the thinning of $W$ obtained by
taking $W$ points and then randomly and independently keeping each of them with
probability $p$. For integers $l \ge 0$ and $0 \le r \le l$ 
let $\pi_{lr}$ denote the binomial probabilities
\begin{align*}
  \pi_{lr}(p)\=\P\bigpar{\Bi(l,p)=r}=\binom lr p^r (1-p)^{l-r}.
\end{align*}
(The understanding here is that $\pi_{00} (p) =1$ for all $p$.)
Thus we have
\begin{equation*}
  \P(W_p=r)=\sum_{l=r}^\infty p_l \pi_{lr}(p).
\end{equation*}
We further define, for given $(p_r)_{r=0}^\infty$, functions 
\begin{align*}
  h(p)&\= 
\E\bigpar{W_p\ett{W_p\ge k}}
=
\sum_{r=k}^\infty \sum_{l=r}^\infty r p_l \pi_{lr}(p),
\\
  h_1(p)&\= \P(W_p\ge k)
=\sum_{r=k}^\infty \sum_{l=r}^\infty  p_l \pi_{lr}(p).
\end{align*}
Note that both $h$ and $h_1$ are increasing in $p$, with
$h(0)=h_1(0)=0$. Note further that $h(1)=\sum_{r=k}^\infty r p_r \le\gl$
and $h_1(1)=\sum_{r=k}^\infty  p_r \le1$, with strict inequalities
unless $p_r=0$ for all $r=1,\dots,k-1$ or $r=0,1,\dots,k-1$, respectively.

The following theorems are our central results, and are key to proving
\refT{T1}. 
See Fernholz and Ramachandran~\cite{fr03,fr04} and in particular
Cooper \cite{c04} for similar results. 

\begin{theorem}
  \label{TD}
Consider the random multigraph \gndx{} for a sequence  $(d_i)_1^n$
satisfying \refCN{C1}.
Let $k\ge2$ be fixed, and let $\crkx$ be the \kcore{} of \gndx.
Let $\pmax$ be the
largest $p\le1$ such that $\gl p^2=h(p)$. 
  \begin{romenumerate}
\item
If $\pmax=0$, \ie{} if $\gl p^2>h(p)$ for all $p\in(0,1]$, then 
$\crkx$ has $\op(n)$ vertices and $\op(n)$  
edges \whp{} (if it exists at all). Furthermore, 
if also $k\ge3$ and $\sumin e^{\ga d_i} =O(n)$ for
some $\ga>0$, then \crkx{} is
empty \whp.
\item
If $\pmax>0$, and further $\gl p^2< h(p)$ for $p$ in some interval
$(\pmax-\eps,\pmax)$, 
then
\whp{} \crkx{} is non-empty, and 
$v(\crkx)/n\pto h_1(\pmax)$, $e(\crkx)/n\pto h(\pmax)/2=\gl \pmax^2/2$.
  \end{romenumerate}
\end{theorem}

\begin{theorem}
  \label{TD2}
If\/ $\sum_i d_i^2=O(n)$ and $\sum_i d_i^3=o(n^{3/2})$
then all the conclusions of \refT{TD} hold also for
the random graph \gnd. 
\end{theorem}

Naturally, the extra condition $\sumin e^{\ga d_i} =O(n)$ in
\refT{TD}(i) implies the extra conditions in \refT{TD2}.

\section{Finding the core}\label{Sfind}

It is well-known (see for instance~\cite{psw96}) that the $k$-core of an
arbitrary finite graph or 
multigraph can be found by removing vertices of degree $\ltx k$, in
arbitrary order, until no such vertices exist. 
It is easily seen that we obtain the same result by removing edges
where one endpoint has degree $\ltx k$, until no such edges remain, and
finally removing all isolated vertices. Again, the order of removal
does not matter, and we will use a randomized choice as follows.

Regard each edge as consisting of two \emph{half-edges}, each
half-edge having one endpoint.
Say that a vertex is \emph{light} if its degree is $\ltx k$, and
\emph{heavy} if its degree is $\gex k$. Similarly, say that a half-edge
is light or heavy when its endpoint is. As long as there is any light
half-edge, choose one such half-edge uniformly at random and remove the edge it
belongs to. 
(Note that this may change the other endpoint from heavy to light, and
thus create new light half-edges.)
When there are no light half-edges left, we stop. 
Then all light vertices are isolated;
the heavy vertices and the remaining edges form
the $k$-core of the original graph.

We apply this algorithm to a random multigraph with given degree
sequence $(d_i)_1^n$. 
Let us observe only the vertex degrees in the resulting multigraph
process, but not the individual edges. In other words, we observe the
half-edges, but not how they are connected into edges.
At each step, we thus select a light half-edge at random. We then
reveal its partner, which is random and uniformly distributed over the
set of all other half-edges. We then remove these two half-edges and
repeat as long as there is any light half-edge. It is clear, by
considering configurations, that this
gives a Markov process (the state at any time $t\ge 0$ is the current
degree sequence); and that at each step, conditioned on the
vertex degrees observed so far, the remaining
multigraph is a random multigraph with the given vertex degrees and
the distribution specified in \refS{Smulti}.

We shall analyse this process of half-edges further in \refS{Sproof}.

\section{Some death processes}\label{S:death}

This section contains some preliminary lemmas that will be used in
our proofs.
We begin with a classical result,
see \eg{} Proposition 4.24 in~\cite{Kallenberg}.
\begin{lemma}[The Glivenko--Cantelli theorem]
  \label{L0}
Let $T_1,\dots,T_n$ be \iid{} random variables with 
distribution function $F(t)\=\P(T_i\le t)$, and let $X_n(t)$ be their
empirical distribution function $\#\set{i\le n: T_i\le t}/n$.
Then $\sup_t|X_n(t)-F(t)|\pto0$ as \ntoo.
\nopf
\end{lemma}

Consider next a pure death process with rate 1; this process starts
with some number of balls whose lifetimes 
are \iid{} rate 1 exponentials $\Exp(1)$.

\begin{lemma}
\label{L1}
Let $N^{(n)}(t)$ be the number of balls alive at time $t$ in a rate $1$
death process with $N^{(n)}(0)=n$.
Then
\begin{equation*}
\sup_{t\ge0} \bigl| N^{(n)}(t)/n-e^{-t}\bigr| \pto0
\qquad \text{as \ntoo}.
\end{equation*}
\end{lemma}

\begin{proof}
$1-N^{(n)}(t)/n$ is the empirical distribution function of the $n$
  lifetimes, which are \iid{}
random variables with the distribution function $1-\et$, $t\ge0$.
Hence the result is an instance of \refL{L0}.
\end{proof}

The death process in \refL{L1} is a Markov process such that, whenever
in state $j$, the process jumps to $j-1$ with intensity $j$, that is
after a random time with distribution $\Exp(1/j)$. We extend this by
allowing the process to take non-integer values as follows.

\begin{lemma}
\label{L2}
Let $\gamma>0$ and $d>0$ be fixed.
Let $N^{(x)}(t)$ be a Markov process such that $N^{(x)}(0) =x$ a.s. and 
transitions are made according to the following rule: whenever
in state $y>0$, the process jumps to $y-d$ with intensity $\gamma y$;
in other words, the waiting time until the next event is $\Exp(1/\gamma y)$ and
each jump is of size $d$ downwards.
Then
\begin{equation*}
\sup_{t\ge0} \bigl| N^{(x)}(t)/x-e^{-\gamma  dt}\bigr| \pto0
\qquad \text{as \xtoo}.
\end{equation*}
\end{lemma}

\begin{proof}
Dividing $N^{(x)}(t)$ by $d$ and $t$ by $\gamma d$ we can rescale the
process, and so we may just as well assume that $d=\gamma=1$. The
process is then the same as the one in 
\refL{L1} if $x=n$ is an integer. In general, consider $\tnxt(t)$, 
a rate 1 death process satisfying
$\tnxt(0) 
= \ceil{x}$. We can couple
$N^{(x)}(t)$ and $\tnxt(t)$ such that both jump whenever the
smaller does, and it is easily seen that under the coupling
$|N^{(x)}(t)-\tnxt(t)|<1$ for all $t$.
The result thus follows from \refL{L1}, which yields
$\sup_{t\ge0} \bigl| \tnxt(t)/\ceil{x}-e^{-\gamma dt}\bigr| \pto0$.
\end{proof}

Now consider $n$ bins with independent rate 1 death processes.
Let $N^{(n)}_j(t)$ denote the number of
balls in bin $j$ at time $t$, where $j=1,\dots,n$ and $t\ge0$. Let
further $U^{(n)}_r(t)\=\#\set{j:N^{(n)}_j(t)=r}$, the number of bins
with exactly $r$ balls, for $r=0,1,\dots$. In what follows we suppress the
superscripts to lighten the notation.

\begin{lemma}\label{LU}
Consider $n$ independent pure death processes $N_i(t)$ with rate $1$
such that $N_i(0)=d_i$, where $(d_i)_1^n$ satisfies \refCN{C1}.
Then, with the above notation, as \ntoo,
\begin{equation*}
\sup_{t\ge0} \sumr 
 r\left|U_r(t)/n-\sum_{l=r}^\infty p_l\pi_{lr}(\et)\right|\pto0.
\end{equation*}
In particular, 
\begin{align}
\sup_{t\ge0}\, \biggl|\sum_{r=k}^\infty r U_r(t)/n-h(\et)\biggr|&\pto0, 
\label{sofie}
\\
\sup_{t\ge0}\, \biggl|\sum_{r=k}^\infty U_r(t)/n-h_1(\et)\biggr|&\pto0.
\label{emma}
\end{align}
\end{lemma}

\begin{proof}
Let $U_{lr}(t)$ be the number of bins that have $l$ balls at time 0
and $r$ balls at time $t$. We shall actually prove the stronger result
\begin{equation}\label{anna}
\sup_{t\ge0} \sum_{l=0}^\infty\sum_{r=0}^l 
 r\left|U_{lr}(t)/n- p_l\pi_{lr}(\et)\right|\pto0.
\end{equation}

First fix integers $l$ and $j $, with $1 \le j\le l$. Consider the 
$u_l\=U_l(0)$ bins that start with $l$ balls. For 
$i=1, \ldots , u_l$ let $\tij$ be the time the $j$-th ball is removed from
the $i$-th such bin.
Then $\#\set{i:\tij \le t}=\sum_{s=0}^{l-j} U_{ls}(t)$.
Moreover, the number of balls remaining in one of these bins at time
$t$ has the distribution $\Bi(l,\et)$, and thus
$\P(\tij\le t) = \sum_{s=0}^{l-j} \pi_{ls}(\et)$. Multiplying by
$u_l/n$ and using \refL{L0}, we obtain that
\begin{equation*}
\sup_{t\ge0} \left|
\frac 1n \sum_{s=0}^{l-j} U_{ls}(t)
- \frac{u_l}n 
\sum_{s=0}^{l-j} \pi_{ls}(\et)\right|\pto0.
\end{equation*}
Further, this convergence trivially holds when $l=0$ or $j=0$.
But $u_l/n\to p_l$ by \refC{C1}(i), and so in fact, for all $j,l\ge0$,
\begin{equation*}
\sup_{t\ge0} \left|
\frac 1n \sum_{s=0}^{l-j} U_{ls}(t)
- p_l
\sum_{s=0}^{l-j} \pi_{ls}(\et)\right|\pto0.
\end{equation*}
Take $j=l-r$ and $j=l-r+1$ and subtract the corresponding
quantities under the absolute value sign to deduce that each term in
\eqref{anna} tends to 0 in probability. 
Hence the same holds for any finite partial sum.

Finally, let $\eps>0$ and let $L$ be such that 
$\sum_L^\infty lp_l<\eps$. By \refCN{C1}(iii), 
$\sum_l lu_l/n\to\gl=\sum_l l p_l$.
Hence also
$\sum_{l\ge L} lu_l/n\to\sum_{l\ge L} l p_l<\eps$.
Consequently, if $n$ is large enough,
$\sum_{l\ge L} lu_l/n<\eps$, and
\begin{align*}
\sup_{t\ge0} \sum_{l=L}^\infty\sum_{r=0}^l 
 r\left|U_{lr}(t)/n- p_l\pi_{lr}(\et)\right|
&\le
\sup_{t\ge0} \sum_{l=L}^\infty\sum_{r=0}^l 
 r\Bigpar{U_{lr}(t)/n+ p_l\pi_{lr}(\et)}
\\&
\le
\sum_{l=L}^\infty
 l\bigpar{u_l/n+ p_l}
<2\eps.
\end{align*}
We conclude that \eqref{anna} holds.
\end{proof}

\section{Proof of \refT{TD}}\label{Sproof}

We continue to analyse the process of vertex degrees in the
core-finding algorithm of \refS{Sfind} applied to a random multigraph
with given degree sequence $(d_i)_1^n$. We regard vertices as bins and
half-edges as balls. The description in \refS{Sfind} thus says
that at each step we remove first one random ball from the set of
balls in light bins
(\ie{} bins with $\ltx k$ balls) and then a random ball without
restriction. We stop when there are no non-empty light bins, and the
$k$-core consists precisely of the heavy bins at the time we stop.

We thus alternately remove a random light ball and a random ball. We
may just as well say that we first remove a random light ball. We then
remove balls in pairs, first a random ball and then a random light
ball, and stop with the random ball leaving no light ball to remove.

We change the description a little by introducing colours. Initially
all balls are white, and we begin again by removing one random light
ball. Subsequently, in each deletion step we first remove a random
white ball and then recolour a random light white ball red; this is repeated
until no more white light balls remain. If we consider only the white
balls, this is evidently the same process as before.

We now run this deletion process in continuous time such that, if there
are $j$ white balls remaining, then 
we wait an exponential time with mean $1/j$ until the next pair of deletions.
In other words, we make deletions at rate $j$. This means that
each white ball is deleted with rate 1 and that, when we delete a
white ball, we also colour a random light white ball red.
Let $L(t)$ and $H(t)$ denote the numbers of light and heavy white
balls at time $t$ respectively; further, let $H_1(t)$ be the
number of heavy bins.

Since red balls are ignored, we may make a final change of rules,
and say that all balls are removed at rate 1 and that, when a
white ball is removed, a random white light ball is coloured red; we
stop when we should recolour a white light ball but there is no
such ball. Note that all heavy balls are white, and that white
balls yield our core-finding process.

Let $\tau$ be the stopping time of this process.
First consider the white balls only.
There are no white light balls left
at $\tau$, so $L(\tau)$ has reached zero. However, let us consider the
last deletion \& recolouring step as completed by redefining
$L(\tau)\=-1$; we then see that $\tau$ is characterized by
$L(\tau)=-1$ and $L(t)\ge0$ for $0\le t<\tau$.
Moreover, the heavy balls left at $\tau$ (which are all white) are
exactly the half-edges in the $k$-core. Hence the number of edges in
the \kcore{} is $\tfrac12H(\tau)$, while the number of vertices is $H_1(\tau)$.

Moreover, if we consider only the total number $L(t)+H(t)$
of white balls in the bins,
ignoring the positions, the process (up to time $\tau$) is
as follows: each ball dies
at rate 1 and upon its death another ball is also sacrificed. 
The process $L(t)+H(t)$ thus is the death process studied in
\refL{L2}, with $\gamma=1$ and $d=2$.
We start with an odd number $2m-1$ of white balls, since we began by
removing one light ball. Consequently, \refL{L2} yields
\begin{equation}\label{erika}
 \sup_{t\le\tau} \bigl|L(t)+H(t)-2m{e^{-2t}}\bigr|
=\op \bigpar{2m} = \op(n).
\end{equation}

Next let us ignore the colours. Our final version of the process
then becomes 
exactly the process studied in \refL{LU},
apart from the initial removal of a light ball which does not affect
the conclusions because, for each $t$, at most two $U_r(t)$ (in the
notation of \refS{S:death}) are
changed by $\pm1$.

Since all heavy balls are
white, we have $H(t)=\sumoo_{r=k} rU_r(t)$ and 
$H_1(t)=\sumoo_{r=k}U_r(t)$.
Hence, 
by \eqref{sofie} and 
\eqref{emma},
\begin{align}
\sup_{t\le\tau} |H(t)/n-h(\et)| & \pto0, \label{jesper}
\\
\sup_{t\le\tau} |H_1(t)/n-h_1(\et)| &\pto0.\label{jesper1}
\end{align}
In particular,
\begin{align}
H(\tau)/n- h(e^{-\tau}) \pto0, 
\quad\text{and}\quad
H_1(\tau)/n- h_1(e^{-\tau}) \pto0.
\label{david}
\end{align}
We deduce from \eqref{erika}, \eqref{jesper} and $2m/n\to\gl$ that
\begin{equation}\label{magnus}
\sup_{t\le\tau} \bigl|L(t)/n+h(\et) -\gl{e^{-2t}}\bigr|
\pto0.
\end{equation}

Assume now that  $t_1$ is a constant independent of $n$ 
with $t_1<-\ln \pmax$.
Then $\pmax <1$ and thus $h(1)<\gl$.
Hence, by continuity, 
$h(p)-\gl p^{2}<0$ on
$(\pmax,1]$, and thus
$h(\et)-\gl e^{-2t}<0$ for $t\le t_1$.
By compactness, 
$h(\et)-\gl e^{-2t}\le -c$ for $t\le t_1$ and some $c>0$.
But $L(\tau)=-1$, so if $\tau\le t_1$ then
 $L(\tau)/n+h(e^{-\tau}) -\gl{e^{-2\tau}}<-c$ and from
\eqref{magnus}
\begin{equation}\label{manne}
  \P(\tau\le t_1)\to0.
\end{equation}

In case (i) we may take any finite $t_1$ here, and hence find
$\tau\pto\infty$. As $h(0)=h_1(0)=0$, \eqref{david} yields that
\begin{align*}
  H(\tau)/n&\pto0, && H_1(\tau)\pto0.
\end{align*}
The first claim now follows, since $v(\crk)=H_1(\tau)$ and
$e(\crk)=H(\tau)/2$. The second claim will follow from \refL{Lempty} below.

In case (ii) we similarly let $t_2\in (-\ln \pmax,-\ln(\pmax-\eps))$.
Then by the hypothesis 
 $h(e^{-t_2}) -\gl{e^{-2t_2}}=c>0$.
If $\tau> t_2$ then $L(t_2)\ge 0$, and thus
 $L(t_2)/n+h(e^{-t_2}) -\gl{e^{-2t_2}}\ge c$. Consequently
\eqref{magnus} implies that
\begin{equation*}
  \P(\tau \ge t_2)\to0.
\end{equation*}
Since we can choose $t_1$ and
$t_2$ arbitrarily close to $-\ln \pmax$, together with \eqref{manne}
this shows that
\begin{equation*}
  \tau\pto - \ln \pmax.
\end{equation*}
Combined with \eqref{david}, this yields 
$H(\tau)/n\pto h(\pmax)$ and
$H_1(\tau)/n\pto h_1(\pmax)$, which proves (ii).
\qed

It remains to prove the following lemma extending 
a result by {\L}uczak \cite{l91}.
\begin{lemma}
  \label{Lempty}
If $k\ge3$ and $\sum_i e^{\ga d_i} =O(n)$, then there exists $\gd>0$
such that
\whp{} \gndx{} has no non-empty \kcore{} with fewer than $\gd n$ vertices.
\end{lemma}

\begin{remark}
The proof below shows the stronger statement that
\whp{} \gndx{} has no non-empty subgraph with fewer than $\gd n$
vertices
and average degree at least $k$.
\end{remark}

We begin the proof of \refL{Lempty} with a sublemma.

\begin{lemma}
  \label{LA}
Consider a set $X$ of $2m$ points and a subset $Y\subseteq X$ with $y$
elements. Let $M$ be a random perfect matching of $X$ and let $Z$ be
the number of pairs in $M$ where both members belong to $Y$.
Then for every real $u\ge0$
\begin{equation}\label{samuel}
  \P(Z \ge u) \le \parfrac{y^2}{mu}^u.
\end{equation}
\end{lemma}

\begin{proof}
Denote the \rhs{} of \eqref{samuel} by $f(u)$. Then either $f(u)\ge 1$
or $f(u)\ge f(\ceil{u})$. Hence it suffices to prove \eqref{samuel}
when $u$ is an integer.
In that case 
\begin{align*}
  \P(Z\ge u)
&
\le 
 \E\binom{Z}{u}=
\binom {y}{2u}\frac{(2u)!}{2^u u!} \frac1{(2m-1)\dotsm(2m-2u+1)}
\\&
=
\binom {y}{2u} \frac{\binom m u} {\binom{2m} {2u}}
\le
\parfrac {y}{2m}^{2u}\binom m u
\le
\parfrac {y}{2m}^{2u}\parfrac{em}{u}^u
=
\parfrac {ey^2}{4mu}^u.
\end{align*}
\end{proof}

\begin{proof}[Proof of \refL{Lempty}]
Let $C$ be such that $\sum_i e^{\ga d_i} \le Cn$.

Consider a set $A$ of $s$ vertices $i_1,\dots,i_s$, and let
$D_A\=\sum_{j=1}^s d_{i_j}$.
If $A$ is the vertex set of the \kcore, it must contain at least
$ks/2$ edges. By \refL{LA}, using the inequality $x\le e^{x}$, the
probability of this event is at most
  \begin{align*}
\parfrac{2 D_A^2}{mks}^{ks/2}   
=
\parfrac{2 ks}{m\ga^2}^{ks/2}\parfrac{\ga D_A}{ks}^{ks}
\le
\parfrac{2 ks}{m\ga^2}^{ks/2}e^{ks \ga D_A/(ks)}.       
  \end{align*}
Summing over all sets $A$ with $s$ vertices, we obtain
\begin{align*}
  \P(v(\crkx)=s)
&\le
\parfrac{2 ks}{m\ga^2}^{ks/2}\sum_{|A|=s} \prod_{i\in A} e^{\ga d_i}
\\&
=
\parfrac{2 ks}{m\ga^2}^{ks/2}\parfrac{n}{s}^s
 \sum_{|A|=s} \prod_{i\in A} \frac sn e^{\ga d_i}
\\&
\le
\parfrac{2 ks}{m\ga^2}^{ks/2}\parfrac{n}{s}^s
  \prod_{i=1}^n \Bigpar{1+\frac sn e^{\ga d_i}}
\\&
\le
\parfrac{2 ks}{m\ga^2}^{ks/2}\parfrac{n}{s}^s
  \exp\Bigpar{\sum_{i=1}^n \frac sn e^{\ga d_i}}
\\&
\le
\parfrac{2 ks}{m\ga^2}^{ks/2}\parfrac{n}{s}^s
  \exp(Cs)
\end{align*}
Since $2m/n\to\gl$, $m>\gl n/3$ for large $n$, so that
\begin{equation}\label{julie}
  \P(v(\crkx)=s)
\le
\Bigpar{
 \parfrac{6 k}{\gl \ga^2}^{k/2}\parfrac{s}{n}^{k/2-1} e^C
}^s.
\end{equation}
Choosing $\gd$ such that 
\begin{align*}
   \parfrac{6 k}{\gl \ga^2}^{k/2}\gd^{k/2-1} e^C =\frac12,
\end{align*}
and considering the cases $s<\ln n$ and $s\ge \ln n$ separately,
it is easily seen that the sum of the \rhs{} of
\eqref{julie} over $s\in [1, \gd n]$ is $o(1)$.
\end{proof}

\section{Proofs of Theorems \refand{TD2}{T1}} \label{Spf2}

\begin{proof}[Proof of \refT{TD2}]
As is well-known, see for instance \cite{b01} and \cite{McKay}, under our
assumptions $\liminf \P(\gndx\text{ is simple})>0$.
Indeed, by considering subsequences we may assume that 
$\sum d_i(d_i-1)/2m\to\mu<\infty$, and then the number of loops and multiple
edges converges, \eg{} by the method of moments, to a
$\Po(\mu/2+\mu^2/4)$ distribution.
Hence the  result follows from \refT{TD} by conditioning on \gndx{}
being simple.
\end{proof}

\begin{proof}[Proof of \refT{T1}]
The degree sequence $(d_i)_1^n$ is now random, but \refCN{C1} holds
for convergence in
probability with $p_r=\P(\Po(\gl)=r)$, see for example \cite[Chapter
III]{b01}. Choosing a suitable coupling of the random graphs \gnl{}
for different $n$, we may thus assume that \refCN{C1} holds a.s.

Further, the vertex degrees $d_i$ all have the same distribution, 
binomial $\Bi(n-1,\gl/n)$ for \gnl{} and hypergeometric for \gnm,
and it follows easily that
$\E \sum_i e^{d_i} = n \E e^{d_1} =O(n)$. This implies that $\sum_i
e^{d_i} = \Op(n)$; by suitable conditioning we may thus assume 
$\sum_i e^{d_i} = O(n)$.
Then \refT{TD2} applies \as{} to \gnl{} or \gnm{} conditioned on the
degree sequence, with
$(p_r)_r=\Po(\gl)$. In the notation of \refS{Smulti},
$W\sim\Po(\gl)$ and so $W_p\sim\Po(\gl p)$;
hence
$h_1(p)=\psik(\gl p)$ and 
\begin{equation*}
  h(p)
=\sum_{j=k}^\infty j \frac{(\gl p)^j}{j!} e^{-\gl p}
=\gl p \psiki(\gl p).
\end{equation*}
Consequently,
\begin{equation*}
  \gl p^2 > h(p) 
\iff
p>\psiki(\gl p)
\iff
\frac{\gl p}{\psiki(\gl p)} > \gl.
\end{equation*}
It then follows that $\pmax=0 \iff \gl<\mu/\psiki(\mu)$ for all $\mu\le\gl$.
Since this inequality holds trivially for $\mu>\gl$, we deduce that
$\pmax=0\iff \gl<\glk$, and so part (i) follows.

Similarly, if $\gl>\gl_k$, $\gl
\pmax=\mu_k(\gl)$, and (ii) follows,
provided we show that $\mu/\psiki(\mu)<\gl$ for $\mu$
slightly less than $\mu_k(\gl)$. This is done in
\refS{S:fixed} below.
\end{proof}

\section{A fixed point equation} \label{S:fixed}

To complete the proof of \refT{T1} we show the following lemma.
Recall that
$\glk\=\min_{\mu>0}\mu/\psiki(\mu)$.

\begin{lemma}
  \label{LP}
  \begin{thmenumerate}
\item
Assume $k\ge3$.
If $\gl>\glk$, then the equation
$\mu/\psiki(\mu)=\gl$ has exactly two positive solutions,
$\muxx-(\gl)$ and $\muxx+(\gl)$, with $0<\muxx-(\gl)<\muxx+(\gl)$;
thus $\mu_k(\gl)=\muxx+(\gl)$.
Moreover,
$\mu/\psiki(\mu)<\gl$ for $\muxx-(\gl)<\mu<\muxx+(\gl)=\mu_k(\gl)$.
\item
Assume $k=2$.
If $\gl>\glk$, then the equation
$\mu/\psiki(\mu)=\gl$ has exactly one positive solution,
$\mu_k(\gl)$, and
$\mu/\psiki(\mu)<\gl$ for $0<\mu<\mu_k(\gl)$.
  \end{thmenumerate}
\end{lemma}

\begin{proof}
  Define $\phi(\mu)\=\psiki(\mu)/\mu$.

For $k=2$, $\psiki(\mu)=1-e^{-\mu}$ so
$\phi(\mu)=\bigpar{1-e^{-\mu}}/\mu$. Hence $\phi$ is (strictly)
decreasing on
$(0,\infty)$ and $\mu/\psiki(\mu)$ is increasing from
$\gl_2=1$ to $\infty$ for $\mu\in(0,\infty)$; the result follows.

For $k\ge3$, the result follows immediately from the lemma below;
note that 
\begin{equation*}
  \glk\=\inf_{\mu>0}\frac1{\phi(\mu)}
  =\frac1{\sup_{\mu>0}\phi(\mu)}.
\end{equation*}
\vskip-\baselineskip
\end{proof}

\begin{lemma}
  \label{LP1}
If $k\ge3$, then $\phi(x)\=\psiki(x)/x$ is unimodal: there is a unique
maximum point $x_0>0$,  $\phi'(x)>0$ for $0<x<x_0$ and $\phi'(x)<0$
for $x>x_0$. Further, $\phi (x)\to 0$ as $x \to 0$ or $x \to \infty$.
\end{lemma}
\begin{proof}
  Note first that $\phi$ is continuously differentiable on
  $(0,\infty)$ with $\phi(x)>0$, and that $\phi(x)=O\bigpar{x^{k-2}}$
  as $x\to0$, and $\phi(x)\le1/x$; hence $\phi(x)\to0$ as $x\to0$ or
  $x\to\infty$. It follows that
$\phi(x)$ attains its maximum at some $x_0>0$.

Also $\psiki'(x)=x^{k-2}e^{-x}/(k-2)!$, and thus
$\psiki(x)/\bigpar{x\psiki'(x)}$ is increasing. Hence
\begin{equation}
  \label{la}
x\frac{d}{dx}\ln\phi(x)
=
x\frac{\psiki'(x)}{\psiki(x)}-1
\qquad
\end{equation}
is decreasing. Since 
$\ln \phi$ attains its maximum at $x_0$, 
$\frac{d}{dx}\ln\phi(x_0)=0$, and it follows from \eqref{la} that 
$\frac{d}{dx}\ln\phi(x)>0$ for $x<x_0$ and
$\frac{d}{dx}\ln\phi(x)<0$ for $x>x_0$.
\end{proof}

\begin{remarks}
The proof shows that $y\mapsto \ln\phi(e^y)$ is strictly concave.

In the
language of discrete dynamical systems, 
see for instance~\cite{d03}, for $k\ge3$,
$\muxx\pm(\gl)$ are the fixed points of 
$f_{\gl}(x)\=\gl\psiki(x)$, and $f_\gl$ undergoes a
saddle-node bifurcation at $\gl=\gl_k$.
\end{remarks}

\section{Further results}

We have studied the $k$-core of a random multigraph with a given
degree sequence. We have determined sufficient conditions on the
asymptotic behaviour of the degree sequence for the $k$-core to be
empty, or at least very small, with high probability. We have also
given sufficient conditions for the multigraph to have a giant
$k$-core such that the sizes of its vertex and edge sets obey a law of
large numbers.

We have further given a new proof that the random graph \gnl{} (and
hence also the 
random graph \gnm) exhibits threshold behaviour. That is, for each
integer $k \ge 3$, there is a value $\glk$ such that, if $\gl < \glk$
then the $k$-core is empty \whp; and if $\gl > \glk$
then the number of vertices and number of edges in the $k$-core are
almost deterministic, and are very large.

We have not discussed the next level of detail. It is possible to
obtain quantitative versions of our results, such as 
large deviation estimates 
and a central limit theorem for the size of the $k$-core. 
Also, one can use our method to study the transition window: how far
above the threshold the edge probability $\gl/n$ must be to ensure
that \gnl{} has a non-empty $k$-core \whp. (Some such results were
already given by Pittel et al.~\cite{psw96}.) These and other issues
will be considered in a forthcoming paper.

Furthermore, it seems possible to adapt the methods of this paper to 
random hypergraphs, but we leave this to the reader.

\newcommand\AAP{\emph{Adv. Appl. Probab.} }
\newcommand\JAP{\emph{J. Appl. Probab.} }
\newcommand\JAMS{\emph{J. \AMS} }
\newcommand\MAMS{\emph{Memoirs \AMS} }
\newcommand\PAMS{\emph{Proc. \AMS} }
\newcommand\TAMS{\emph{Trans. \AMS} }
\newcommand\AnnMS{\emph{Ann. Math. Statist.} }
\newcommand\AnnPr{\emph{Ann. Probab.} }
\newcommand\CPC{\emph{Combin. Probab. Comput.} }
\newcommand\JMAA{\emph{J. Math. Anal. Appl.} }
\newcommand\RSA{\emph{Random Struct. Alg.} }
\newcommand\SPA{\emph{Stoch. Proc. Appl.} }
\newcommand\ZW{\emph{Z. Wahrsch. Verw. Gebiete} }
\newcommand\DMTCS{\jour{Discr. Math. Theor. Comput. Sci.} }

\newcommand\AMS{Amer. Math. Soc.}
\newcommand\Springer{Springer}
\newcommand\Wiley{Wiley}

\newcommand\vol{\textbf}
\newcommand\jour{\emph}
\newcommand\book{\emph}
\newcommand\inbook{\emph}
\def\no#1#2,{\unskip#2, no. #1,} 

\newcommand\webcite[1]{\hfil\penalty0\texttt{\def~{\~{}}#1}\hfill\hfill}
\newcommand\webcitesvante{\webcite{http://www.math.uu.se/\~{}svante/papers/}}

\def\nobibitem#1\par{}

\end{document}